\documentclass{article}
\usepackage{url}
\usepackage{amsthm}
\usepackage{amssymb}
\usepackage{enumerate}
\usepackage{lmodern}
\usepackage{natbib}
\usepackage[english]{babel}

\usepackage{amsthm}
\usepackage{graphicx}
\usepackage{amsmath}
\usepackage{amsthm}
\usepackage{amssymb}
\usepackage{algorithmic}
\usepackage[ruled,vlined]{algorithm2e}
\usepackage{braket}
\usepackage{comment}
\usepackage{mathtools}
\usepackage[utf8]{inputenc}
\usepackage[mathscr]{eucal}
\usepackage{fancyhdr}
\usepackage{fancybox}
\usepackage{xcolor}
\usepackage{longtable}
\usepackage{dsfont}
\usepackage[colorlinks=true, allcolors=blue]{hyperref}
 {
     \theoremstyle{plain}
     
 }
 
\newtheorem {Proposition}{Proposition}[section]
\newtheorem {Lemma}[Proposition] {Lemma}
\newtheorem {Theorem}[Proposition]{Theorem}
\newtheorem {Corollary}[Proposition]{Corollary}

\newtheorem {Definition}{Definition}[section]

\usepackage{graphicx}
\usepackage{amsfonts}
\usepackage{amsmath}
\usepackage[utf8]{inputenc}
\usepackage{url}
\usepackage{color}
\usepackage{amsthm}
\usepackage{amssymb}
\usepackage{enumerate}
 \usepackage{paralist}
\usepackage{graphicx}

\usepackage{mathtools}
\usepackage{color}
\RequirePackage{amsthm,amsmath,amsfonts,amssymb}

\def\N{\mathbb{N}}

\def\R{\mathbb{R}}

\numberwithin{equation}{section}

\newcommand{\ssubset}{\subset\joinrel\subset}

\usepackage{amssymb}

\begin{document}



\begin{center}
{\Large 
Regularity of center-outward distribution functions in non-convex domains} \\
\vspace{0.5cm} 
Eustasio del Barrio$^{1}$\\
Alberto González Sanz$^{1,2}$
\\
\vspace{0.2cm} 
$^1$  IMUVa, University of Valladolid \\
$^2$ IMT,  Université Paul Sabatier
\end{center}

\begin{abstract}
For a probability $\rm P$ in $\R^d$ its center outward distribution function $\mathbf{F}_{\pm}$, introduced in \cite{Hallin_17} and \cite{Hallin_annals}, is a new and successful concept of multivariate distribution function based on mass transportation theory.
This work  proves, for a probability $\rm P$ with density locally bounded  away from zero and infinity in its support, the continuity of the center-outward map on the interior of the support of $\rm P$ and  the continuity of its inverse, the quantile,  $\mathbf{Q}_{\pm}$. This relaxes  the convexity assumption in \cite{BaGoHa}. Some important  consequences of this continuity are Glivenko-Cantelli type theorems and characterisation of weak convergence by the stability of  the center-outward map. \end{abstract}





\section{Introduction\label{sec:1}}

There is a rich literature defining generalizations of the well-known univariate distribution function. This is not only one the most powerful tools from a
probabilistic point of view (it completely determines a probability measure or characterizes weak convergence of probability measures, for instance). Additionally, and, perhaps, more importantly, distribution functions and their empirical versions underlie an efficient approach to nonparametric and semiparametric statistical inference, namely, the approach based on ranks. We refer to  \cite{Hallin2006SEMIPARAMETRICALLYER} for a detailed account on these aspects. The success of rank based statistical methods was, until recently, constrained to univariate data.

Given this considerations it cannot be a surprise the large amount of literature devoted to different extensions of the univariate distribution functions to higher dimension and the related concept of ranks. These generalizations come from different motivations related to depth, elliptic models or componentwise extensions among other (\cite{Halin2015,Hallin2010MultivariateQA,HallinReview}). However, since the univariate distribution function is deeply related with the natural ordering of $\R$, all of them have had to deal with the same problem, the lack of canonical order in dimension higher than one.
From the point of view of statistical inference, none of the extensions before \cite{Hallin_17}  and \cite{Hallin_annals} was able to produce a
meaningful version of mutivariate ranks. These last two works introduced a multivariate extension of the distribution function which, in turn, allows to construct semiparametrically efficient statistical procedures for multivariate data \cite{HallinReview}. The idea for this extension is based on the mass transportation characterization of the the distribution function $F$---it is the optimal transport map from the associated probability  measure $P$ to the uniform distribution on $[0,1]$, or, more generally, the unique gradient of a convex map that pushes forward $\rm P$ to the uniform on $[0,1]$. With a view to a multivariate extension, in which the left-to-right order is meaningless, it makes sense to replace the uniform on $[0,1]$ by the uniform on the unit ball centered at the origin. This generalizes to higher dimension: the center-outward distribution function introduced in \cite{Hallin_annals}  is unique gradient of a convex map that pushes forward a probability $\rm P$ to the spherical uniform measure on the unit ball $ \mathrm{U}_d$, which means an independent random choice in the sphere (the direction) and the radius. We will write $\mathbf{F}_{\pm}$ for this center outward distribution function (we refer to subsection~\ref{sec:notat2} below for a formal definition). The choice of the spherical uniform as the reference measure is not completely inocent. It is the only one which defines proper ranks and signs, which yield the desired statistical properties---useful for many applications such as \cite{HallinVARMA,HallinKendall,HallinRankTest,quantileReg}, among others. On the other hand, this choice poses also a technical challenge when trying to guarantee its existence as a proper, single-valued map.

From a mathematical point of view, a first difficulty with the center outward distribution function comes from the fact that, even under smoothness assumptions on $\rm P$, $\mathbf{F}_{\pm}$ is not uniquely defined at every point. For a probability $\rm P\in \mathcal{P}(\R^d)$ with a density with respect the $d$-dimensional Lebesgue measure $\ell_d$, a celebrated result of \cite{Maccan} guarantees that there exists a map which is, almost everywhere, the gradient of a convex map pushing forward $\rm P$ to $ \mathrm{U}_d$. The map is $\rm P$- almost surely unique.
As it is well known, a convex funtion is almost everywhere differentiable and points of nondifferentiability are those in which the subgradient is not a sigleton. This means that one could circunvent the lack of definition at some points by thinking of $\mathbf{F}_{\pm}$ as a possibly multivalued map. A significant recent contribution in this line is \cite{Segers}, which analyzes the connection between weak convergence of probability measures and convergence of center-outward distribution functions as multivalued maps in the topology of Painlevé-Kuratowski.

From a statistical point of view, it is a remarkable fact that the univariate distribution function can be uniformly approximated by its empirical version, as the Glivenko-Cantelli Theorem states. Uniform convergence, even locally, cannot hold at points in which the center-outward distribution function is multivalued. Hence, it is important to provide (a) sufficiently general conditions under which the center-outward distribution function is single valued and (b) uniform approximation guarantees for the empirical versions. These two tasks are the main goal of this paper.

Since for gradients of convex functions continuity and single-valuedness are equivalent properties (see e.g. Theorem~25.1 in \cite{Rockafellar} or  Theorem 12.25 in \cite{RockWets98}), task (a) can be reformulated in terms of sufficiently general conditions under which the center-outward distribution function is continuous. Hence, the problem is related to regularity conditions for monotone measure preserving maps, which can be explored using the powerful regularity theory of Caffarelli (cf. \cite{Caffarelli1990ALP,Caffarelli1991SomeRP,Ca3}) for solutions of the Monge-Ampère equation. However, the classic theory applies to probabilities with densities bounded away from $0$ and infinity in convex supports. Unfortunately, the spherical uniform $ \mathrm{U}_d$ that led to the meaningful multivariate extension has a density unbounded around the origin. Hence, neither the classic results, nor the most general one of \cite{Erausquin} apply to our setting. \\

Predating this paper, two works provide regularity results for the center-outward distribution function; \cite{FigalliCenter} for probabilities with densities locally bounded away from $0$ and $\infty$ in all of $\R^d$ and  \cite{BaGoHa} for probabilities supported in a convex set and with densities locally bounded away from $0$ and $\infty$ in their supports. The result of  \cite{BaGoHa} generalizes the one in \cite{FigalliCenter}. In both cases not only the regularity of $\mathbf{F}_{\pm}$ is given, but also that of the quantile map, $\mathbf{Q}_{\pm}$, its inverse: both maps are homeomorphisms between the support of $\rm P$ minus the median set, $\mathbf{Q}_{\pm}({\bf 0})$, and the punctured unit ball $\mathbb{B}_d-\{\mathbf{0}\}$. The median set $\mathbf{Q}_{\pm}({\bf 0})$ is a compact convex negligible set, but possibly not a singleton, which means that the continuity of $\mathbf{Q}_{\pm}$ requires additional assumptions. But, as we prove below, continuity of the center-outward distribution holds under less restrictive hypotheses on the probability  $\rm P$ and, in fact, we consider the following local boundedness condition on the density:

\noindent \emph{Assumption \textbf{A}}:\emph{ $\rm P$  is absolutely continuous with respect to $\ell_d$, with density $p$. Denote $ \mathcal{X}:=\operatorname{int}\left(\operatorname{supp}(\rm P) \right).$ For every $R>0$,  there exist constants $0<\lambda_R\leq  \Lambda_R$ such that
\begin{equation}\label{upperlowwer}
 \lambda_R  \leq p({\bf x}) \leq
 \Lambda_R \quad\text{ for all ${\bf x}\in{\mathcal X}\cap R\,\mathbb{B}_d$}.\bigskip
\end{equation}}
Under this assumption, the main result in this paper, Theorem~\ref{TheoremMain}, shows that:
 \begin{enumerate}
     \item The center-outward distribution map ${\mathbf F}_\pm$ is continuous in $\mathcal{X}$ and  the center-outward quantile map ${\mathbf Q}_\pm$ is continuous in the set ${\mathbf F}_\pm(\mathcal{X})-\{\mathbf{0}\}$. Moreover, $\mathbf{Q}_\pm$  and $\mathbf{F}_\pm$  are homeomorphisms between  ${\mathbf F}_\pm(\mathcal{X})-\{\mathbf{0}\}$ and ${\mathcal X}-\mathbf{Q}_\pm(\mathbf{0})$, inverses of each other.
     \item If $p\in \mathcal{C}^{k,\alpha}_{\operatorname{loc}}(\mathcal{X})$, for some $k\in \N$ and $\alpha\in (0,1)$, then $\mathbf{Q}_\pm$  and $\mathbf{F}_\pm$  are diffeomorphisms of class $\mathcal{C}^{k,\alpha}_{\operatorname{loc}}$  between  ${\mathbf F}_\pm(\mathcal{X})-\{\mathbf{0}\}$ and ${\mathcal X}-\mathbf{Q}_\pm(\mathbf{0})$.
     \item  If  $\mathcal{X}$ is convex, then  ${\mathbf F}_\pm(\mathcal{X})=\mathbb{B}_d$.
 \end{enumerate}
The last conclusion is the main result in \cite{BaGoHa}, but Theorem~\ref{TheoremMain} generalizes it to a wider class of probability measures not necessary supported on a convex set.\\

In the usual statistical data analysis setup $\rm P$ is not directly observed, but there is access to an i.i.d.  sample $\mathbf{X}_1, \dots,\mathbf{X}_n$ with underlying law $\rm P$. Since the empirical measure $\rm {\rm P_n}$ does not have a density, its center-outward distribution function cannot be defined directly in terms of optimal transport maps. To circunvent this issue \cite{Hallin_annals} proposes the use of a discretization of the spherical measure and defines the empirical distribution function through an interpolation procedure. The consistency of this estimator,
$\hat{\mathbf F}_{n,\pm}$, ---the Glivenko–Cantelli type theorem---relies basically on the regularity of ${\mathbf F}_\pm$. Theorem~\ref{TheoremMain} has as a consequence Corollary~\ref{colo:Gliv}, which yields uniform consistency of the empirical center-outward distribution function in compact sets $K$ contained in its support, i.e.
 \begin{equation*}
     \sup_{\mathbf{x}\in K}| \hat{\mathbf F}_{n,\pm}(\mathbf{x})-{\mathbf F}_\pm(\mathbf{x})|\longrightarrow 0, \ \ a.s.
 \end{equation*}
 
We would like to remark at this point that our approach to Corollary~\ref{colo:Gliv} comes from a more general consequence of  Theorem~\ref{TheoremMain} which has probabilistic implications of interest by itself. Suppose that $\mu_{m}$, with continuous center-outward distribution function ${\mathbf F}_\pm^m$, converges weakly to $\rm P$, with center-outward distribution function ${\mathbf F}_\pm$ that we also assume that is continuous. Then (see Theorem~\ref{theorem:Weakconv} below), for every compact subset $K$ of the interior of the suppport of $\rm P$
\begin{equation}\label{Glivenco}
     \sup_{\mathbf{x}\in K}| {\mathbf F}_\pm^m(\mathbf{x})-{\mathbf F}_\pm(\mathbf{x})|\longrightarrow 0.
 \end{equation}
This means that the center-outward maps characterize weak convergence of probability measures with continuous center-outward distribution.
\\

The rest of the paper is organized as follows. Section \ref{sec:notat} sets the notation used throughout the paper. Section \ref{sec:mainresult} presents the main result, Theorem~\ref{TheoremMain}, with a discussion about its assumptions and improvements with respect to previous works. For ease of readability its proof is postponed to Section~\ref{sec;proof}. Section \ref{sec:conseq}  illustrates the main consequences  of Theorem~\ref{TheoremMain} in two settings, one related to statistical consistency and the other on with a characterization of weak convergence.

\subsection{Formal definition of the center-outward distribution function.}\label{sec:notat}

The formal definition of the center-outward distribution function relies, ultimately, on some usual concepts from convex analysis. We refer to \cite{Rockafellar} or \cite{RockWets98} for further details, but introduce here some required elements from that theory. Let $f$ be a lower semicontinuous convex function and $\mathbf{x}\in \operatorname{dom}(f)$. The subdifferential of $f$ in $\textbf{x}$ is defined as
$$
\partial f(\textbf{x}):=\{\textbf{y}\in R^d: f(\textbf{z})\geq f(\textbf{x})+\left\langle \textbf{y},\textbf{z}-\textbf{x} \right\rangle, \ \text{for all $\textbf{z}\in \R^d$} \}.
$$
The Legendre transform of $f$ is the convex function
$
f^*(\mathbf{y}):=\sup_{\mathbf{x}\in \R^d} (\langle\mathbf{x},\mathbf{y} \rangle-f(\mathbf{x})),\quad \mathbf{y}\in \mathbb{R}^d.
$

The spherical uniform distribution $\mathrm{U}_d$ over the open $d-$dimensional unit ball, $\mathbb{B}_d$, is the probability measure with density $ u_d(\mathbf{x})=\frac 1 {a_d|\mathbf{x}|^{d-1}} I\big[\mathbf{x}\in \mathbb{B}_d-\{\mathbf{0}\}\big],$ where $a_d=2\pi^{d/2}/\Gamma(d/2)$. This corresponds to a random (uniform) choice of a direction on the   shpere and, independently, a random (also uniform) choice of a distance to the origin. We say that a measurable map $T$ pushes forward $\rm P\in \mathcal{P}(\R^d)$ into $\mathrm Q\in \mathcal{P}(\R^d)$ if $\rm{P}(A)=\rm{Q}(T^{-1}(A))$, for every Borel measurable set $A\subset \R^d$.
It follows from \cite{Maccan}, that there exists a convex function, $\psi$, such that $\nabla \psi$ pushes $\mathrm{U}_d$ forward to $\rm P$. We call $\mathbf{Q}_\pm:=\nabla \psi$ the center-outward quantile map associated to $\mathrm P$ and note that is is well defined at almost every point in the open unit ball. $\psi$ is uniquely defined (up to an additive constant, that we fix by setting $\psi(\mathbf{0})=0$) in the open unit ball.
We proceed to extend $\psi$ to all of $\mathbb{R}^d$ by setting
\begin{align*}
    \text{$\psi(\mathbf{u}):=\liminf_{\mathbf{z}\to\mathbf{u},|\psi(\mathbf{z})|<1} \psi(\mathbf{z})$ if~$|\mathbf{u}|=1$ and~$\psi(\mathbf{u}):=+\infty$ for $\mathbf{u}\notin \bar{\mathbb{B}}_d$.}
\end{align*}
Now we set $\varphi:=\psi^*$, the convex conjugate of $\psi$, namely,
\begin{equation}\label{def_phi}
\varphi(\mathbf{x})=\sup_{|u|<1} (\langle\mathbf{x},\mathbf{u} \rangle-\psi(\mathbf{u})),\quad \mathbf{x}\in \mathbb{R}^d.
\end{equation}
We observe that \textit{(i)} the domain of $\varphi$ is $\mathbb{R}^d$, \textit{(ii)} $\varphi$ is 1-Lipschitz (being the sup of a family of 1-Lipschitz functions) and
\textit{(iii)} $\nabla \varphi$ pushes $\rm P$ forward to $\mathrm{U}_d$. The center-outward distribution function is defined now as $\mathbf{F}_{\pm}=\nabla \varphi$. This is a well-defined, single valued function at almost every $\mathbf{x}\in\mathbb{R}^d$. In fact, setting $\mathbf{F}_{\pm}=\partial \varphi$, we would have, with the above construction, a valid definition of the, possibly multivalued, center-outward distribution function over all of $\mathbb{R}^d$, also beyond the support of $\rm P$. We note also that the fact that $\varphi$ is 1-Lipschitz implies that
\begin{equation}\label{contencion1}
\partial \varphi(\mathbb{R}^d)\subset \bar{\mathbb{B}}_d.
\end{equation}
For the sake of illustration, we discuss two examples of center-outward distribution functions. First, we take $\mathrm{P}=\mathrm{U}_d$. In this case $\psi(\mathbf{u})=\frac{|u|^2}2$, $|u|\leq 1$ ($\nabla\psi (\mathbf{u})= \mathbf{u}$ in $\mathbb{B}_d$), $\psi(\mathbf{u})=+\infty$, $|\mathbf{u}|>1$. A simple computation yields $\varphi(\mathbf{x})=
\frac{|\mathbf{x}|^2} 2$, $|\mathbf{x}|\leq 1$; $\varphi(\mathbf{x})= |\mathbf{x}|- \frac 1 2$, $|\mathbf{x}|> 1$. We conclude that the center-outward distribution function of the spherical uniform distribution equals
$$\mathbf{F}_\pm(\mathbf{x})=\left\{
\begin{matrix}
\mathbf{x}, & & |\mathbf{x}|\leq 1,\\
\frac{\mathbf{x}}{|\mathbf{x}|}, & & |\mathbf{x}|> 1.
\end{matrix}
\right. $$
For our second example we take as $P$ the law of $\mathbf{X}=\mathbf{U}+\text{sgn}(U_1)\mathbf{e}_1$, where $\mathbf{U}$ is a random vector with spherical uniform distribution, $\mathrm{U}_d$, $U_1$ denotes the first component of $\mathbf{U}$ and $\mathbf{e}_1$ is the first element of the canonical basis of $\mathbb{R}^d$. Now $\psi(\mathbf{u})=\frac{|\mathbf{u}|^2}2+|u_1|$, $|\mathbf{u}|\leq 1$, where we write $\mathbf{u}=(u_1,\ldots,u_d)$. It is convenient to write $\mathbf{x}=(x_1,\mathbf{x}_2)$, with $\mathbf{x}_2\in\mathbb{R}^{d-1}$. Routine computations show that $\varphi(\mathbf{x})=\frac {|\mathbf{x}-\text{sgn}(x_1)\mathbf{e}_1|^2}2$ if $|x_1|\geq 1, |\mathbf{x}-\text{sgn}(x_1)\mathbf{e}_1|\leq 1$;
$\varphi(\mathbf{x})=|\mathbf{x}-\text{sgn}(x_1)\mathbf{e}_1| -\frac 1 2$ if $|x_1|\geq 1, |\mathbf{x}-\text{sgn}(x_1)\mathbf{e}_1|> 1$;
$\varphi(\mathbf{x})=\frac{|\mathbf{x}_2|^2} 2$ if $|x_1|\leq 1, |\mathbf{x}_2|\leq 1$ and, finally, $\varphi(\mathbf{x})=|\mathbf{x}_2|-\frac 1 2$ if $|x_1|\leq 1, |\mathbf{x}_2|> 1$. We conclude that
$$\mathbf{F}_\pm(\mathbf{x})=\left\{
\begin{matrix}
\mathbf{x}-\text{sgn}(x_1)\mathbf{e}_1, & & \text{if }|x_1|\geq 1, |\mathbf{x}-\text{sgn}(x_1)\mathbf{e}_1|\leq 1,\\
\frac{\mathbf{x}-\text{sgn}(x_1)\mathbf{e}_1}{|\mathbf{x}-\text{sgn}(x_1)\mathbf{e}_1|}, & & \text{if }|x_1|\geq 1, |\mathbf{x}-\text{sgn}(x_1)\mathbf{e}_1|>1,\\
(0,\mathbf{x}_2), & & \text{if }|x_1|\leq 1, |\mathbf{x}_2|\leq 1,\\
\frac{(0,\mathbf{x}_2)}{|\mathbf{x}_2|}, & & \text{if }|x_1|\leq 1, |\mathbf{x}_2|> 1.
\end{matrix}
\right. $$

\subsection{Notation.}\label{sec:notat2}

We end this Section with a few lines on notation.
The absolute value of real numbers and the Euclidean norm of vectors are both denoted by $|\cdot |$. Set $\mathbf{x}_0\in \R^d$ and $\tau>0$, the notation $\mathbf{x}_0+\tau\mathbb{B}_d=\{\mathbf{x}\in \R^d: \ |\mathbf{x}-\mathbf{x}_0|<\tau \} $ refers to the open Euclidean ball of $\R^d$ centered in $\mathbf{x}_0$ and radius $\tau$.  For subsets $A, B$, of a metric space we write $A\ssubset B$ to denote that $A$ is compactly contained in $B$, that is, there exists a compact set $K$, contained in the interior of $B$, such that $A\subset K$.
 
The space of Borel probabilities on $\R^d$ is denoted by $\mathcal{P}(\R^d)$. The support of a probability measure is given by $\operatorname{supp}(\cdot)$. The Lebesgue measure in $\R^d$ is denoted by $\ell_d$. Let $\mu,\nu $ be two Borel measures, the relation $\mu\ll\nu$ means that $\mu$ is absolutely continuous with respect to  $\nu$.  For two sets $A,B\,\subseteq \mathbb{R}^d$, its Hausdorff distance (see \cite{RockWets98}) is defined as $\mathrm{d}_H(A,B):= \max\big\{\sup_{{\bf a}\in A}\inf_{_{\scriptstyle{\bf b}\in B}}\vert {\bf a}-{\bf b}\vert ,\,\sup_{{\bf b}\in B} \inf_{_{\scriptstyle{\bf a}\in A}}\vert {\bf a}-{\bf b}\vert \big\}.$ The closed convex hull of a set $A$ is denoted by $ \bar{\operatorname{co}}(A)$. Let $A$ be an open set of $\R^d$,  the notation $ \mathcal{C}^{k,\alpha}_{\operatorname{loc}}(A)$, for $k\in \N$ and $\alpha\in (0,1)$, refers to the space of Hölder continuous functions

\section{Continuity of center-outward distribution functions.}\label{sec:mainresult}

As noted in the Introduction, \cite{FigalliCenter} and \cite{BaGoHa} provide conditions under which the center-outward distribution function is continuous. They only cover convexly supported probabilities, but their goal, however, was not only to prove continuity of  $\mathbf{F}_{\pm}$, but also that of its inverse---the quantile function $\mathbf{Q}_{\pm}$. Here we show that the assumption of convexity can be relaxed, at least to preserve the continuity of ${\mathbf F}_\pm$. In fact, when the probability is not convexly supported, we cannot guarantee the continuity of $\mathbf{Q}_{\pm}$. 
The next result guarantees continuity of $\mathbf{F}_\pm$ inside the interior of the support of $\rm P$, $\mathcal{X}$, under \emph{Assumption \textbf{A}} (see \eqref{upperlowwer}). Regularity outside the interior of the support requires some additional smoothness at the boundary $\partial\mathcal{X}$. We cannot expect $\mathbf{F}_{\pm}$ to be continuous without making further assumptions at least on the shape of the support of $\rm P$. To our knowledge, all the existing literature dealing with this  boundary regularity assumes at least the convexity of $\mathcal{X}$ (cf. \cite{CaffarelliBound1,CaffarelliBound2} or, more generally, see Remark~4.25 in \cite{Fi} and references therein).
Continuity of the quantile holds on the set ${\mathbf F}_\pm(\mathcal{X})-\{\mathbf{0}\}$. With these assumptions the regularity of the center-outward map and its quantile inverse is described in the following theorem, which is the main contribution of this work.
\begin{Theorem}\label{TheoremMain}
Let  $\rm P\in \mathcal{P}(\R^d)$ be a probability with density $p$ satisfying \emph{Assumption \textbf{A}}. Then the following assertions hold.
\begin{enumerate}
\item The center-outward distribution map ${\mathbf F}_\pm$ is continuous in $\mathcal{X}$ and  the center-outward quantile map ${\mathbf Q}_\pm$ is continuous in the set ${\mathbf F}_\pm(\mathcal{X})-\{\mathbf{0}\}$. Moreover, $\mathbf{Q}_\pm$  and $\mathbf{F}_\pm$  are homeomorphisms between  ${\mathbf F}_\pm(\mathcal{X})-\{\mathbf{0}\}$ and ${\mathcal X}-\{\mathbf{Q}_\pm(\mathbf{0})\}$, inverses of each other.
\item If $p\in \mathcal{C}^{k,\alpha}_{\operatorname{loc}}(\mathcal{X})$, for some $k\in \N$ and $\alpha\in (0,1)$, then $\mathbf{Q}_\pm$  and $\mathbf{F}_\pm$  are diffeomorphisms of class $\mathcal{C}^{k,\alpha}_{\operatorname{loc}}$  between  ${\mathbf F}_\pm(\mathcal{X})-\{\mathbf{0}\}$ and ${\mathcal X}-\{\mathbf{Q}_\pm(\mathbf{0})\}$.
\item  If  $\mathcal{X}$ is convex, then  ${\mathbf F}_\pm(\mathcal{X})=\mathbb{B}_d$.
\end{enumerate}
\end{Theorem}

To prove Theorem~\ref{TheoremMain}, we  need to avoid the singularity of the spherical distribution. To do it we define an auxiliary restricted Monge-Amp\`ere measure (also called Brenier solutions of the Monge-Amp\`ere equation in \cite{FiKim}). We provide   Alexandrov  upper bounds (see Lemma \ref{lem:Alexandrov}). The lower control is provided in Lemmas \ref{lem:bounded2} and \ref{lem:previo1}. Finally Lemma~\ref{lem:convex} proves the injectivity of the quantile map which yields the continuity of $\mathbf{F}_{\pm}$. The continuity of the quantile $\mathbf{Q}_{\pm}$---in a set which does not contain the singular point---follows the arguments of  \cite{FiKim}. 

There exist numerous works dealing with the regularity of the optimal transport map, see \cite{Caffarelli1990ALP,Caffarelli1991SomeRP,Caffarelli1992,FiKim,Fi}. In most cases the assumptions include some, at least local, lower and upper bound on the density of the measures, which does not hold for the spherical uniform measure. We refer to \cite{BaGoHa} for further discussion on the peculiarities of this case. Recently, \cite{FigalliDoubling} extended the work of \cite{Erausquin} to the setup of \textit{doubling measures}.  A measure $\mu$ is said to be locally doubling (\cite{Jhaveri2022}) if for every ball ${\mathbb B}_d$, there is a constant $C\geq 1$, such that
\begin{equation}
\label{doubling}
\mu(S)\leq C \mu({\textstyle\frac{1}{2}}S),
\end{equation}
for any convex set $S\subset {\mathbb B}_d$  with center (of mass) in $\operatorname{supp}(\mu)$. Here $\frac{1}{2}S$ denotes  the dilation of $S$  with respect to its center by $\frac{1}{2}$. However, as we show below (Proposition \ref{doubling_prop}), the spherical uniform measure is not doubling for $d>2$ and, thereforem, Theorem~\ref{TheoremMain} is not covered by any exsisting result in the literature.
Furthermore, we must emphasize on the fact that, to our knowledge, this is the first continuity result for the center-outward distribution function without convexity assumptions on the shape of the support of the associated probability measure.

\begin{Proposition}\label{doubling_prop}
The spherical uniform ${\rm U}_d$ is not locally doubling for $d>2$.
\end{Proposition}
\begin{proof}
Let $S_r=[-\frac{1}{12},\frac{7}{12}]\times[-r,r]^{d-1}$ be a rectangle with center of mass $(\frac{1}{4}, 0,\dots, 0)$. Note that $S_r\subset\mathbb{B}_d$, for $r$ small enough. On the one hand, since $\frac{1}{2}S_r={[\frac{1}{12},\frac{5}{12}]\times[-\frac{r}{2},\frac{r}{2}]^{d-1}}$, lies at a positive distance, { $\frac{1}{12}$,} from $\bf 0$, there exists $C_1>0$ such that
$$ {\rm U}_d({\textstyle \frac{1}{2}}S_r){ \leq} C_1 \ell_d({\textstyle \frac{1}{2}}S_r)={\textstyle\frac{C_1}{3\cdot 2^{d-1}}r^{d-1}}.$$
On the other hand,  $r\mathbb{B}_d\subset S_r-\frac{1}{2}S_r$, for small $r$ , and
$ {\rm U}_d({ r}\mathbb{B}_d)=r.$ Therefore, if ${\rm U}_d$ were doubling (note that since ${\rm U}_d$ is bounded, the term `locally' is redundant),
$${\textstyle\frac{C_1}{3\cdot 2^{d-1}}r^{d-1}}\geq {\textstyle \frac{1}{C} }{\rm U}_d (S_r)\geq  {\textstyle \frac{1}{C} } {\rm U}_d({ r} \mathbb{B}_d)= {\textstyle \frac{r}{C}}$$
for small enough $r$, which yields a contradiction if $d>2$.
\end{proof}

To end this section we remark that continuity of $\mathbf{F}_\pm$ holds over all of $\mathbb{R}^d$ under the assumption of convex $\mathcal{X}$ (Proposition 2.3 in \cite{Hallin_annals}). Without convexity, proving continuity of $\mathbf{F}_\pm$ in $\mathbb{R}^d$ would imply continuity at the boundary of the support of $P$, $\partial\mathcal{X}$. This is a delicate issue, for which the existing literature is very scarce. We should cite here the remarkable exception of \cite{MIURA2021107603} and the references therein for a broader view on the topic.

\section{Consequences}\label{sec:conseq}

\subsection{Glivenko–Cantelli theorem}
This section shows that, under \emph{Assumption \textbf{A}}, the center-outward map can be consistently  estimated from the sample. Corollary~\ref{colo:Gliv} is the analogous to the one dimensional  Glivenko–Cantelli theorem for the distribution function. \\ 

Let $\mathbf{X}_1, \dots,\mathbf{X}_n$  be a sample of i.i.d. random variables with law $\rm P\in \mathcal{P}(\R^d)$. Denote as ${\rm P_n}$ its empirical measure. The definition of the empirical version of the center outward map requires a \textit{discretized approximation} of $\mathrm{U}_d$, call it ${\rm U}_d^n$,  which satisfies ${\rm U}_d^n\xrightarrow{w}\mathrm{U}_d$.
\begin{Definition}
We say that a continuous map $\hat{\mathbf F}_\pm:\R^d\rightarrow\mathbb{B}_d$ is an empirical regular interpolation of the center outward map if there exists a convex function $\hat{\varphi}:\R^d\rightarrow\R$ such that $\hat{\mathbf F}_\pm=\nabla \hat{\varphi}$,
$\hat{\mathbf F}_\pm\sharp {\rm P_n}={\rm U}_d^n$ and $\partial\hat{\varphi}(\R^d)\subset\overline{\mathbb{B}}_d$.
\end{Definition}
Not every discretized approximation of $\mathrm{U}_d$ yields the existence of empirical regular interpolations. We refer to \cite{Hallin_annals}  for one which guarantees this existence. For any of them the following result is consequence of Theorem~\ref{TheoremMain}.
\begin{Corollary}\label{colo:Gliv}
Let  $\rm P\in \mathcal{P}(\R^d)$ be  such that \emph{Assumption \textbf{A}} holds, then
 \begin{equation}\label{Glivenco3}
\sup_{\mathbf{x}\in K}| {\mathbf F}_\pm^n(\mathbf{x})-{\mathbf F}_\pm(\mathbf{x})|\longrightarrow 0, \ \ a.s.
 \end{equation}
 for every compact $K\subset \mathcal{X}$, and 
  \begin{equation}\label{quantile}
     \sup_{\mathbf{u}\in M}| {\mathbf Q}_\pm^n(\mathbf{u})-{\mathbf Q}_\pm(\mathbf{u})|\longrightarrow 0, \ \ a.s.
 \end{equation}
 for every compact $M\subset \mathbf{F}_{\pm}(\mathcal{X})-\{\mathbf{0}\}$.
 Moreover, if $\mathcal{X}$ is convex, then 
 \begin{equation}\label{Glivenco2}
     \sup_{\mathbf{x}\in \R^d}| {\mathbf F}_\pm^n(\mathbf{x})-{\mathbf F}_\pm(\mathbf{x})|\longrightarrow 0, \ \ a.s.
 \end{equation}

\end{Corollary}
The proof of the second part of Corollary~\ref{colo:Gliv}  is a direct consequence of Theorem~\ref{theorem:Weakconv}. The first part can be derived form the argument to prove uniform convergence of the maps on compact sets within the proof of Theorem~\ref{theorem:Weakconv}.

 \subsection{Characterization of weak convergence.}
The characterization of weak convergence---in the sense of $\mu_m\xrightarrow{w}\rm P$ iff $F^m(x)\rightarrow F(x)$ for every continuity point of $F$---is one of the main properties of the univariate distribution function. Moreover, it is also well known that if $F$ is  everywhere continuous this convergence is uniform (see for instance Exercise~14.8 in \cite{billing}). This, in fact, is translated verbatim to our multivariate case. The techniques of the proof are inspired by those in Proposition~3.3 in \cite{Hallin_annals} and Theorem~2.8. in \cite{BaLo} (see also \cite{Segers}), up to some special technical details. Suppose that $\mu_{m}$ has  continuous center-outward distribution function ${\mathbf F}_\pm^m$ and   converges weakly to $\rm P$. We assume also that ${\mathbf F}_\pm$ is continuous.
\begin{Theorem}\label{theorem:Weakconv}
Suppose that $\mu_{m}\in \mathcal{P}(\R^d)$ has a continuous center-outward distribution function ${\mathbf F}_\pm^m$, for all  $m\in \N$, and   converges weakly to $\rm P\in \mathcal{P}(\R^d)$. Assume also that ${\mathbf F}_\pm$ is continuous. Then
 \begin{equation*}
     \sup_{\mathbf{x}\in \R^d}| {\mathbf F}_\pm^m(\mathbf{x})-{\mathbf F}_\pm(\mathbf{x})|\longrightarrow 0.
 \end{equation*}
\end{Theorem}
In particular, Theorem \ref{theorem:Weakconv} holds for any probability satisfying \emph{Assumption \textbf{A}} with convex support. For those who satisfy only  \emph{Assumption \textbf{A}} the convergence is no longer uniform on $\R^d$ but uniform on the compact sets of the interior of the support of $\rm P$.
\begin{proof} We write $\mathbf{F}_{\pm}^m$ for the center-outward distribution map, defined as the extension of the (continuous) gradient of $\varphi_m$.\\
\textit{(a) Weak convergence of the coupling}:
Relative compactness (for the weak convergence) of the sequence $\{ \mu_m \}_m$ entails (see \cite{Villani2003}) that of the class of probabilities on $\R^d\times \R^d$ with first marginal belonging to $\{ \mu_m \}_n$ and the second equal to $\mathrm{U}_d$. The measures defined by the relations
$
    \pi_m=\left(Id\times{\mathbf F}_\pm^m\right)\sharp\mu_m
$, for $m\in \N$,
belong to such a relatively compact class. As a consequence, it has at least one limit point. Let $\pi\in \mathcal{P}(\R^d\times \R^d)$ be one of them. Lemma~9 in \cite{Maccan} (ii) implies that the marginals of $\pi$  are ${\rm P}$ and $\mathrm{U}_d$. Moreover, since the support of $\pi_m$ is contained in the graph of ${\mathbf F}_\pm^m$ (in the subdifferential of $\varphi$), then it is cyclically monotone  (see \cite{Rocka1966}). Lemma~9 in \cite{Maccan} (i) implies that also $\pi$ is supported on a cyclically monotone set.  Corollary 14 in \cite{Maccan} yields---since ${\rm P}$  and $\mathrm{U}_d$ are uniformly continuous w.r.t. Lebesgue measure---the existence of an unique measure cyclically monotone supported and with marginals  ${\rm P}$ and $\mathrm{U}_d$. As a consequence, we have the limit $\pi_n\xrightarrow{w}\pi$.  \\
\textit{(b) Convergence of the maps on the compact sets}:
Convergence of the couplings implies, using Theorem~2.8 in \cite{BaLo}, the convergence
\begin{align}
    \label{conv:pot}
    \varphi_m\longrightarrow\varphi, \ \ \ \mathrm{P}-a.e.
\end{align}
Note that, Theorem~2.8 in \cite{BaLo} only yields the convergence \eqref{conv:pot} under the assumption of convergence in Wasserstein $2$ topology, but this is only used to show the convergence $\pi_n\xrightarrow{w}\pi$. Now we claim that \eqref{conv:pot} entails uniform convergence on the compact subsets of $\mathbb{B}_d$. To prove this let $K\subset\mathbb{B}_d$ be a compact set and note that the extension described in Section~\ref{sec:notat} guarantees that  $\nabla\varphi_m(\R^d)=\mathbf{F}_{\pm}^m(\R^d) \subset \bar{\mathbb{B}}_d.$ In consequence the modulus of continuity is uniformly controlled and the sequence $\{\varphi_m \}_m$ is equicontinuous. Finally, Arzelà–Ascoli theorem proves the claim.

Now we can prove convergence of $\mathbf{F}_{\pm}^m$ on the compact sets of $\R^d$. Let $K\subset\R^d$ be a compact set, Theorem~25.7 in \cite{Rockafellar},  yields
 \begin{equation*}
     \sup_{\mathbf{x}\in K}| {\mathbf F}_\pm^m(\mathbf{x})-{\mathbf F}_\pm(\mathbf{x})|\longrightarrow 0.
 \end{equation*}
\textit{(c) Uniform convergence on $\R^d$}: In this part we will extract subsequences of subsequences at some point. For the ease of reading we will keep the same notation for these subsequences. Assume that convergence is not uniform.  This implies the existence of $\epsilon>0$ and a subsequence $\{\mathbf{x}_{m}\}_m$ such that
 \begin{equation}\label{contradic}
     | {\mathbf F}_\pm^m(\mathbf{x}_m)-{\mathbf F}_\pm(\mathbf{x}_m)\rangle|>\epsilon, \ \text{for all}\ m\in \N.
 \end{equation}
  and $| \mathbf{x}_{m}|\rightarrow\infty$.
The sequence  $\{\frac{\mathbf{x}_{m}}{|\mathbf{x}_{m}|}\}_m$ belongs to the compact sphere, therefore, along subsequences, $\frac{\mathbf{x}_{m}}{|\mathbf{x}_{m}|}\rightarrow \mathbf{v}\in \mathbb{S}_{d-1}$.
Since $ {\mathbf F}_\pm^m(\mathbf{x}_m), {\mathbf F}_\pm(\mathbf{x}_m)\in \bar{\mathbb{B}}_d$, for all $m\in \N$, they have, respectively, two limit points $\mathbf{y}_1,\mathbf{y}_2 \in \bar{\mathbb{B}}_d$.
 Now we show that $\mathbf{y}_1=\mathbf{v}$, the same can be done, even more easily, with $\mathbf{y}_2$. This means $\mathbf{y}_1=\mathbf{y}_2$, which contradicts \eqref{contradic}. Set $\lambda>0$ and note that the monotony of the subdifferential implies that, for $n$ big enough, the bound
$$\left\langle {\mathbf F}_\pm^m(\mathbf{x}_m)- {\mathbf F}_\pm^m(\lambda\frac{\mathbf{x}_{m}}{|\mathbf{x}_{m}|}),\frac{\mathbf{x}_{m}}{|\mathbf{x}_{m}|}  \right\rangle \geq 0. $$
Taking limits and using the convergence of the maps on the compact sets, we deduce that
\begin{equation}
    \label{laststep}
    \left\langle \mathbf{y}_1- {\mathbf F}_\pm(\lambda\mathbf{v}),\mathbf{v}  \right\rangle \geq 0.
\end{equation}
Proposition 3.1 in \cite{BaGoHa}  shows that the limit of any sequence  $\mathbf{y}_n\in \partial\varphi(|\mathbf{x}_{m}|\mathbf{v})$ is $\mathbf{v}$. Therefore, taking limits for $\lambda\rightarrow+\infty$ in \eqref{laststep}, we obtain that $  \left\langle \mathbf{y}_1- \mathbf{v} ,\mathbf{v}  \right\rangle \geq 0,$ which implies that $ \left\langle \mathbf{y}_1,\mathbf{v}  \right\rangle \geq |\mathbf{v} |^2$. Since $|\mathbf{y}_1|=1$, then $\mathbf{y}_1=\mathbf{v}$.

\end{proof}

\section{Proof of Theorem~\ref{TheoremMain}}\label{sec;proof}
Our proof relies on a careful use of some classical tools from convex analysis and regularity theory for solutions of Monge-Ampère equations.
For those non familiarized with this topic, we refer to \cite{Fi} and references therein. Here we summarize some useful concepts in the course of the proof:
\begin{itemize}
\item A convex set $\Omega\subset \R^d$ is said to be normalized if $\mathbb{B}_d\subseteq \Omega \subseteq d\,\mathbb{B}_d$. For each open bounded convex set $\Omega$ there exists a unique invertible affine transformation $L$   normalizing $\Omega$ (see Lemma A.13 in \cite{Fi}). We refer to $L$ as the {\it normalizing map} and to   $L(\Omega)$ as the {\it normalized version} of $\Omega$.
\item Let $g$ a convex function. The section of $g$ centered at ${\bf x}$ with slope $\mathbf{p}$ and height $t$ is the set
\begin{align}
S(\mathbf{x},\mathbf{p},t):=\{ \mathbf{z}\in \R^d:\ g(\mathbf{z})< l_{\mathbf{x},\mathbf{p},t}(\mathbf{z})\},
\end{align}
where $l_{\mathbf{x},\mathbf{p},t}(\mathbf{z}):=g(\mathbf{x})+\left\langle\mathbf{p},\mathbf{z}-\mathbf{x}\right\rangle+t$.
\item The Monge-Amp\`ere measure associated with a convex function, $\varphi$, is  defined as
$$\mu_\varphi(E):=\ell_d\Big(\partial \varphi(E)\Big)$$
for every Borel set $E{\subseteq} \mathbb{R}^d$. It can be checked that $\mu_\varphi$ is indeed a locally finite Borel measure on $\cal X$.
\end{itemize}
{ We use the notation $C(a_1,\dots,a_m)$ to refer to some positive constant depending only on $a_1,\dots,a_m$. }
\\

We recall from subsection \ref{sec:notat} that $\nabla \varphi$ pushes $\rm P$ forward to ${\rm U}_d$. This implies that $\nabla\varphi(\mathbf{x})\in \mathbb{B}_d$ for almost every $\mathbf{x}\in\mathcal{X}$ and, since  $\mathbb{B}_d$ is convex, we see that \begin{align}\label{cont_ball}
\partial\varphi(\mathcal{X})\subseteq \bar{\mathbb{B}}_d.
\end{align}
A similar conclusion would hold for $\psi$ if $\mathcal{X}$ were convex. Without this assumption we can only guarantee
\begin{align}\label{contained_convex_hull}
\partial\psi(\mathbb{B}_d)\subseteq \bar{co}\left(\mathcal{X}\right),
\end{align}
where $\bar{co}\left(\mathcal{X}\right)$ denotes the closed convex hull of $\mathcal{X}$.\\

The following result, contained in \cite{BaGoHa}, yields the main properties of the Monge-Amp\`ere measure associated with $\varphi$.
\begin{Lemma}\label{MongeAmpereMeasures}
Let  $\rm P$ be a probability measure with  density $p$ supported on  $\mathcal{X}\subseteq \mathbb{R}^d$. Then $\mu_\varphi$ is absolutely continuous with respect to $\ell_d$ and, for every Borel $A\subseteq \mathbb{R}^d$,
$$\mu_{ \varphi } (A)={a_d} \int_{A}p({\bf x})|\nabla \varphi( {\bf x})|^{d-1}d{\bf x}.$$
\end{Lemma}
Lemma \ref{MongeAmpereMeasures} shows that $\mu_{ \varphi }$ has a density w.r.t. $\ell_d$. If $\mathcal{X}$  is not convex,  this is no longer guaranteed for $\mu_{ \psi }$. Instead of working with the whole Monge-Amp\`ere measure associated to $\psi$, we can use its restriction to $\mathcal{X}$, namely,
\begin{align}
\mu_{ \psi }^\mathcal{X}(A):=\ell_d(\partial \psi(A)\cap \mathcal{X}).
\end{align}
This restriction is related with the so-called \emph{Brenier solution} of Monge-Amp\`ere equation in \cite{FiKim}.
We observe that the next result holds with no assumptions on the shape of $\mathcal{X}$.
\begin{Lemma}\label{Restr_MongeAmpereMeasures}
Let  $\rm P$ be a probability measure on $\mathbb{R}^d$ such that \emph{Assumption \textbf{A}} holds. Then $\mu_\psi^\mathcal{X}$ is absolutely continuous with respect to $\ell_d$ and, for every Borel set $B\subseteq \mathbb{B}_d$,
$$\mu_{ \psi }^\mathcal{X}(B)=\frac{1}{a_d}\int_{B}\frac{1}{p(\nabla \psi({\bf y}))|{\bf y}|^{d-1}}d{\bf y}.$$
\end{Lemma} 
\begin{proof}
Let $B\ssubset \R^d$ be a Lebesgue negligible set. We assume, without loss of generality that $\psi(B)$ is bounded. Then, for some $R>0$,
\begin{align*}\label{eq:abs_cont2}
\mu_{ \psi }^{\mathcal{X}} (B)=\ell_d(\partial \psi(B)\cap \mathcal{X}) &\leq \frac{1}{\lambda_R}\int_{\partial \psi(B)}p({\bf x})d{\bf x}
= \frac{1}{\lambda_R}\int_{B}u_d({\bf y})d{\bf y}=0.
\end{align*}
Hence $\mu_{ \psi }^{\mathcal{X}}\ll\ell_d$. Now, using  Theorem 4.8 in \cite{Villani2003}, for every Borel set $B\subset \mathbb{B}_d$ we have,
\begin{equation}\label{eq:abs_cont_dens2}
\mu_{ \psi }^{\mathcal{X}}(B)=\int_{B} \frac{u_d({\bf y})}{p(\nabla \psi( {\bf y}))}d{\bf y}=\frac 1{a_d} \int_{B}\frac{1}{p(\nabla \psi({\bf y}))|{\bf y}|^{d-1}}d{\bf y}. 
\end{equation}
\end{proof}
The distribution $\mathrm{U}_d$ cannot be treated with the standard Caffarelli regularity theory because of its singulary at the origin. The next result ensures that away from it the restricted Monge-Amp\`ere measure is locally bounded away from $0$ and infinity.

\begin{Lemma}\label{lem:bounded2}
If $\rm P$ satisfies \emph{Assumption \textbf{A}} and $M$ is a compact subset of~$\mathbb{B}_d$ such that $\mathbf{0}\not \in M$, then, for every Borel set $A\subseteq M$, we have
\begin{align*}
	\alpha_M\ell_d(A)\leq \mu_{\psi}^{\mathcal{X}}(A)\leq \Lambda_M\ell_d(A).
\end{align*}
\end{Lemma}
\begin{proof}
Note that since $M$ is compact, also its image through $\partial\psi$ (see Lemma A.22 in \cite{Fi}). Then there exists $R>0$ such that $\partial\psi(M)\subset R\mathbb{B}_d$.  \emph{Assumption \textbf{A}} implies the existence of some $\alpha_M,\Lambda_M>0$ such that
$\alpha_M\leq p(\mathbf{q})\leq \Lambda_M$ for all $\mathbf{q}\in \partial\psi(M)\cap \mathcal{X}$. 
In view of Lemma~\ref{Restr_MongeAmpereMeasures}, we have 
$\mu_{ \psi }^\mathcal{X}(A)=\frac{1}{a_d}\int_{A}\frac{1}{p(\nabla \psi({\bf y}))|{\bf y}|^{d-1}}d{\bf y},$ and, since $d_H(\mathbf{0},M)=\delta>0$, then
$$\frac{1}{a_d \Lambda_M}\ell_d(A)\leq \mu_{ \psi }^\mathcal{X}(A)\leq \frac{1}{a_d \alpha_M\delta}\ell_d(A), \ \text{for every Borel set $A\subseteq M$.}$$
\end{proof}
One of the key ideas for the regularity of a convex function $f$ associated with a Monge-Amp\`ere measure such that $$\lambda\ell_d\ll\mu_f\ll\Lambda \ell_d, \quad \text{with $\lambda, \Lambda>0,$}$$  is the invariance of the bounds for the normalized version of $f$ with respect to a section $S= \{f\leq \ell\}\subseteq \mathbb{B}_d$,
$$\mathbf{z}\mapsto| \operatorname{det}(L)|^{\frac{2}{d}}(f-\ell)(L^{-1}{{\bf z}}),$$
meaning that its associated Monge-Amp\`ere  measure is also upper and lower bounded by $\Lambda$ and $\lambda$.
The following result shows that in our setup the lower bound of  Lemma \ref{lem:bounded2}  remains invariant. The upper bound, however, cannot be extended to the normalized solution as in the standard theory.
\begin{Lemma}\label{lem:previo1}
Let $S= \{\psi\leq \ell\}\subseteq \mathbb{B}_d$ be a section of $\psi$ such that $\partial \psi(S)\subset R\,\mathbb{B}_d$, and $L$ be the map which normalizes $S$. Then there exists $\lambda_R>0$ such that the convex function $$v_{L}({{\bf z}}):=| \operatorname{det}(L)|^{\frac{2}{d}}(\psi-\ell)(L^{-1}{{\bf z}})$$ satisfies $\ell_d(\partial v_L(A))\geq \lambda_R\ell_d(A)$ for every $A\subset L(S)$.
\end{Lemma}
\begin{proof}
Let ${\bf z}$ be such that $L^{-1}{{\bf z}}$ is a continuity  point of $\partial\psi$, then we can compute 
$$\text{det}(\operatorname{D}^2(v_L))({\bf z})=\text{det}(\operatorname{D}^2(\psi))(L^{-1}{\bf z})\geq \lambda_R.$$
Even though $\mu_{ \psi }$ might not have a density, Lemma \ref{Restr_MongeAmpereMeasures} ensures  that $\mu_{ \psi }^{\mathcal{X}}$ has, namely,  $\mu_{ \psi }^{\mathcal{X}}(A)=\int_{A} f({\bf x})d{\bf x}$ for some $f$. This implies that if $\Omega:=\{{\bf z}\in \mathbb{B}_d:\  \partial\psi({\bf z})\subseteq\mathcal{X} \}$ then
$$
           \ell_d(\partial v_L(A))= \int_A f(L^{-1}{\bf x})d{\bf x}, \ \ \text{for all}\ \ A\subseteq L( \Omega).
   $$
   Hence, for each $A\subseteq L(S) $, we can split $A=(A\cap L\Omega )\cup(A\cap L\Omega^c)$ and 
   $$\ell_d(\partial v_L(A))\geq \ell_d(\partial v_L((A\cap L\Omega )))=\int_{A\cap L\Omega } f(L^{-1}{\bf x})d{\bf x}\geq \lambda_R \ell_d(A\cap L\Omega ).$$
   Since the set of differentiable points of $\psi$  is of full (Lebesgue) measure in $\mathbb{B}_d$, the same happens for $\Omega$ and the proof follows.
\end{proof}
The strategy now is to show that $\psi$ is strictly convex in the set $\mathbf{F}_{\pm}(\mathcal{X})$, which implies that $\varphi$ (its convex conjugate) is differentiable in $\mathcal{X}$.
\begin{Lemma}\label{lem:convex}  Let  $\rm P$ be a probability measure satisfying \emph{Assumption \textbf{A}}. Then ${\psi}$ is  strictly convex on
\begin{align*}
\Theta := \{ {\bf z}\in \mathbb{B}_d: \ \partial \psi ({\bf z})\cap \mathcal{X}\neq \emptyset\}.
\end{align*}
\end{Lemma}
\begin{proof}
Assume, on the contrary,  that there exists ${{\bf y}}\in \mathbb{B}_d$ and
${\bf t}\in\partial \psi({{\bf y}})\cap \mathcal{X}$ such that, setting $\ell({{\bf z}}):={\psi}({{\bf y}})+\left\langle {{\bf t}}, {{\bf z}}-{{\bf y}}\right\rangle$, the convex   set~$\Sigma:=\{{{\bf z}}:\,  {\psi}({{\bf z}})=\ell({{\bf z}})\}$ is not a singleton. By subtracting an affine function, we can assume $\mathbf{t}=\mathbf{0}\in\mathcal{X} $, ${\psi}({{\bf y}})~\!=~\!0$ and~${\psi}({{\bf z}})\geq 0$ for all ${{\bf z}}$; then, {
$\Sigma=\{{{\bf z}}:\, {\psi}({{\bf z}})=0\} = \{{{\bf z}}:\, {\psi}({{\bf z}})\leq 0\}$}, which  is closed since ${\psi}$ is lower semicontinuous. Also, by adding the convex function $w({{\bf z}}):=\frac 1 2 (|{{\bf z}}|-1)_+^2$ (note that~${\psi}={\psi}+w$ on $\bar{\mathbb{B}}_d$), 
we can assume that $\Sigma\subset \bar{\mathbb{B}}_d$. Being compact and convex,  $\Sigma$ equals the closed convex hull of its extreme points; as a consequence,  it 
 must have at least two exposed points (otherwise it would be empty or a singleton). {Let $\bar{{{\bf y}}}\in\bar{\mathbb{B}}_d\setminus\{\mathbf{0}\}$ be one of them. }  The following result proves that $\bar{{{\bf y}}}$ cannot belong to the sphere. 
 \begin{Lemma}\label{lem:boundary}
Let $\rm P$ be a probability measure satisfying
\emph{Assumption \textbf{A}}. Then $(\partial \varphi )(\mathcal{X})\cap \mathbb{S}_{d-1}=\emptyset$.
\end{Lemma}
\begin{proof}
Note that the convexity assumption in the proof of Lemma 2.2 in \cite{BaGoHa} is redundant.
\end{proof}
As a consequence of Lemma \ref{lem:boundary} $\Sigma\cap \mathbb{S}_{d-1}=\emptyset$ and thus $\bar{{{\bf y}}}\in \mathbb{B}_d\setminus\{\mathbf{0}\}$.
Without any loss of generality, let us assume  that $\bar{{\bf y}}=a\textbf{e}_1$, where $\textbf{e}_1$ stands for the first vector in the canonical basis of $\mathbb{R}^d$ and $a\in (0,1)$.  Now, let $\Pi=\{\mathbf{z}:\ \langle\mathbf{u},\mathbf{z}-\bar{\mathbf{y}}\rangle=0\}$, with $|\mathbf{u} |=1$, be the supporting plane of $\Sigma$ passing through $\bar{\mathbf{y}}$, where we assume that $ \Sigma\subset \{\mathbf{z}:\ \langle\mathbf{u},\mathbf{z}-\bar{\mathbf{y}}\rangle\geq 0\}$.  Let $\mathbf{c}$ be the projection of $\mathbf{0}$ onto the line $\{a\textbf{e}_1+t\, \mathbf{u}:\ t\in \R \}$. Then, either  $\mathbf{c}$ is in the same halfspace of $\Sigma$ and $c_u:=\langle\mathbf{u},\mathbf{c}-\bar{\mathbf{y}}\rangle\geq 0$, or it belongs to its complement. Both cases can be handled simultaneously by defining
\begin{align*}
c_u^{+}=\left\lbrace\begin{array}{cc}
    c_u & \  \text{if}\ c_u>0 \\
    \frac{1}{2} &\ \text{if}\ c_u\leq 0
\end{array} \right., \ \psi_{\epsilon}({{\bf z}}):={\psi}({{\bf z}})-\epsilon (-\langle\mathbf{u},\mathbf{z}-\bar{\mathbf{y}}\rangle+\frac{1}{2}c_u^{+})
\end{align*}
and $S_{\epsilon}:=\{{{\bf z}}:  \psi_{\epsilon}({{\bf z}}) < 0\}.$
Now
\begin{align}\label{Hconv}
S_{\epsilon}\longrightarrow \Sigma \cap \big\lbrace {{\bf z}}\in \R^d:\ \langle\mathbf{u},\mathbf{z}-\bar{\mathbf{y}}\rangle\leq \frac{1}{2}c_u^{+}\big\rbrace
\end{align}
in Hausdorff distance. We point out that (see Figure~\ref{fig:explain}) $$\mathbf{0}\not \in  \Sigma \cap \big\lbrace {{\bf z}}\in \R^d:\ \langle\mathbf{u},\mathbf{z}-\bar{\mathbf{y}}\rangle\leq \frac{1}{2}c_u^+\big\rbrace.$$
\begin{figure}[h!]
    \centering
\includegraphics[width=8cm,height=6cm]{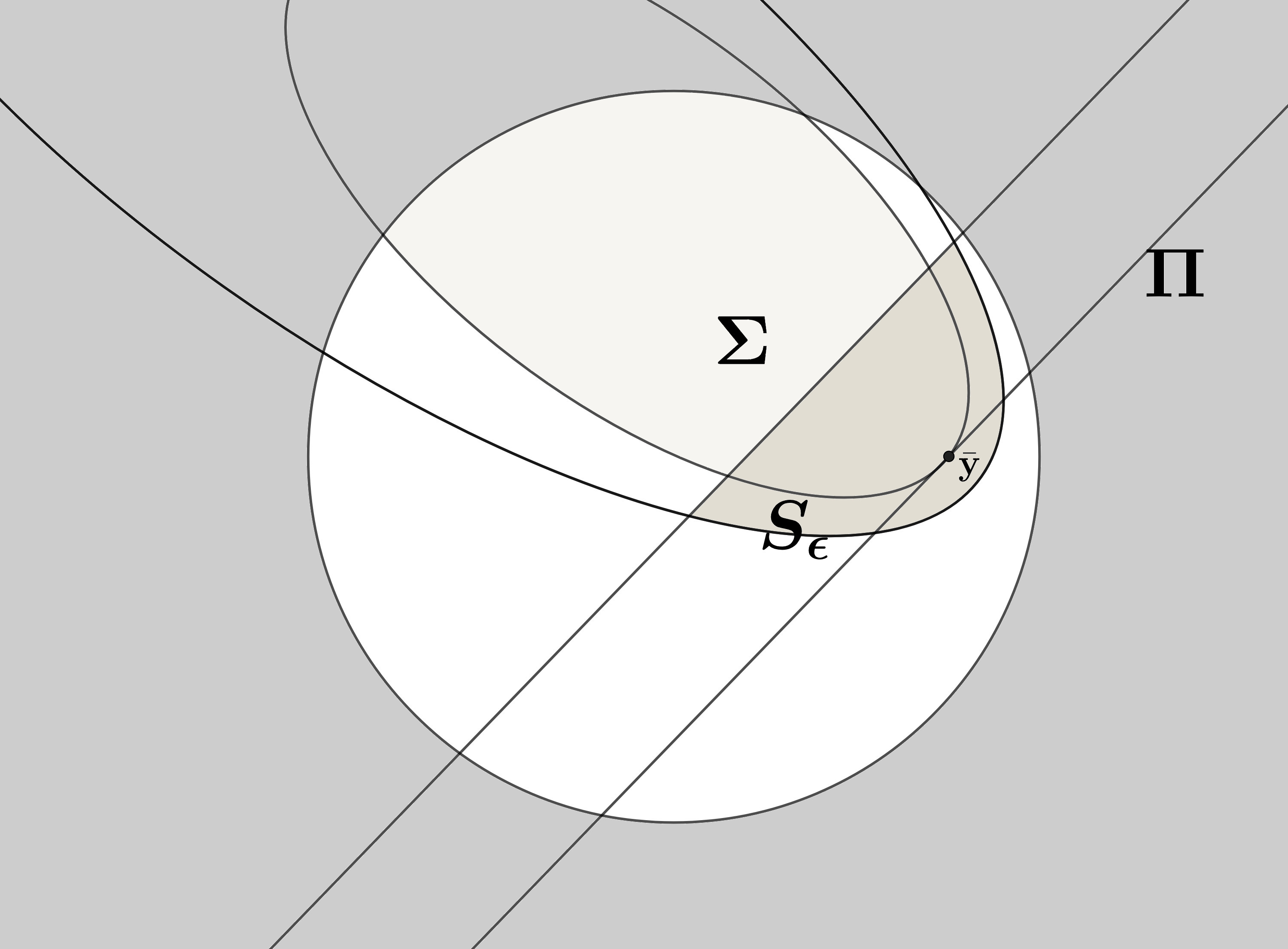}
    \caption{Graphical summary of the argument in the proof of Lemma \ref{lem:convex}. In order to be far away from the singular point (the origin), we project it to the line perpendicular to $\Pi$ and passing through $\bar{\mathbf{x}}$ and we define $\psi_{\epsilon}$ in a way that is $0$ outside of a set $S_{\epsilon}$ which does not touches the origin. This enables to use standard arguments.}
    \label{fig:explain}
\end{figure}

This implies that there exists  $\epsilon_0>0$ such that $ S_{\epsilon}\cap \delta\mathbb{B}_d =\emptyset$, for every ${\epsilon}<\epsilon_0$ and some $\delta>0$.  Moreover, since $\Sigma\cap \mathbb{S}_{d-1}=\emptyset$, there exists $\rho>0$ such that $d_H\left(S_{\epsilon},\mathbb{S}_{d-1}\right)>\rho$, for ${\epsilon}<\epsilon_0$. This separation from the boundary yields (see Lemma A.22 in \cite{Fi})  the relatively compactness of $\partial\psi_{\epsilon}(S_{\epsilon})\subset R\, \mathbb{B}_d$. In other words, there exists $R>0$ such that
\begin{equation}
    \label{compactness}
    \text{$\partial\psi_{\epsilon}(S_{\epsilon})\subset R\, \mathbb{B}_d$, for all ${\epsilon}<\epsilon_0$.}
\end{equation}
Consider the normalizing map~$L_{\epsilon}$ of $S_{\epsilon}$; denote by $v_{\epsilon}({\bf z}):=| \operatorname{det}(L_{\epsilon})|^{\frac{2}{d}}(\psi_{\epsilon}-\frac{a}{2})(L_{\epsilon}^{-1}{\bf z})$. We observe that, 
\begin{align*}
   \lim_{\epsilon\rightarrow 0} \frac{\psi_{\epsilon}(\bar{{\bf y}} )}{\min_{S_{\epsilon}}\psi_{\epsilon}({{\bf y})}}&=\lim_{\epsilon\rightarrow 0}\frac{\psi(\bar{{\bf y}})-\epsilon (\frac{1}{2}c_u^{+}) }{\min_{S_{\epsilon}}\psi({{\bf y})}+\epsilon (\langle\mathbf{u},\mathbf{z}-\bar{\mathbf{y}}\rangle-\frac{1}{2}c_u^{+})}\\
  &=\lim_{\epsilon\rightarrow 0}\frac{\epsilon \frac{1}{2}c_u^{+} }{\max_{S_{\epsilon}}-\psi({{\bf y})}-\epsilon (\langle\mathbf{u},\mathbf{z}-\bar{\mathbf{y}}\rangle-\frac{1}{2}c_u^{+}))}\\
    &\geq \lim_{\epsilon\rightarrow 0}\frac{ \frac{1}{2}c_u^{+}}{\max_{S_{\epsilon}} (\langle\mathbf{u},\mathbf{z}-\bar{\mathbf{y}}\rangle+\frac{1}{2}c_u^{+}))}\geq 1\, .
\end{align*} 
Then, the fact that linear transformations preserve the ratio of the distances between parallel hyperplanes  (see  \cite{Fi} p.78) and  the previous limit imply that
\begin{align}\label{eq:contrastFi}
 \lim_{\epsilon\rightarrow 0}\mathrm{d}_H\left(L_{\epsilon} \bar{\mathbf{y}}, \partial( S_{\epsilon}^L) \right)= 0\quad\text{and}\quad   \lim_{\epsilon\rightarrow 0} \frac{\psi_{\epsilon}(\bar{{\bf y}} )}{\min_{S_{\epsilon}}\psi_{\epsilon}({{\bf y})}} = 1.
\end{align}

Now we will upper bound $|\psi_{\epsilon}(\bar{{\bf y}} )|$ by  $ \mathrm{d}_H\left(L_{\epsilon} \bar{\mathbf{y}}, \partial( S_{\epsilon}^L)\right)$. The next result, based on the well known  Alexandrov upper bound, completes the task.
\begin{Lemma}\label{lem:Alexandrov} 
For each ${\bf y} \in S_{\epsilon}$ such that $\partial\psi ({\bf y})\cap \{{\bf x}: d_H({\bf x}, \partial \mathcal X)>\delta \}\neq \emptyset$, we have that
\begin{align*}
|\psi_{\epsilon}({\bf y})|^d \leq   C(d,\delta,R)\ell_d\left(\partial \psi (S_{\epsilon})\cap \mathcal{X}\cap R\,\mathbb{B}_{d}\right)\ell_d(S_{\epsilon})d_H(L_{\epsilon}(\mathbf{y}), \partial L_{\epsilon}(S_{\epsilon})),
\end{align*}
\end{Lemma}
\begin{proof}
Set ${\bf x}\in \partial\psi ({\bf y})\cap \{{\bf x}: d_H({\bf x}, \partial X)>\delta \}$, then there exist $\delta'\leq\delta$ such that ${\bf x}+\delta'\,\mathbb{B}_{d}\subset \mathcal{X}\cap R\,\mathbb{B}_{d}$. Since $\delta'$ depends on $R$ and $\delta$, we are not losing generality if we assume $\delta=\delta'$.
We consider the convex function $C_{{\bf y}}: S_{\epsilon}\rightarrow \R$ defined as the unique satisfying $C_{{\bf y}}({\bf y})=\psi_{\epsilon}({\bf y})$, $C_{{\bf y}}({\bf z})=\bf 0$ for ${\bf z} \in \partial S_{\epsilon}$, and such that ${\bf z}\mapsto C_{{\bf y}}({\bf y}-{\bf z})$ is 1-homogeneous ($C_{{\bf y}}({\bf y}-s{\bf z})=s C_{{\bf y}}({\bf y}-{\bf z})$, for all $s\in \R-\{0 \}$). Moreover, each supporting plane of $C_{{\bf y}}$ in ${\bf y}$ can be moved since it touched $\psi$ from above inside $S_{\epsilon}$, then it is true that $ \partial C_{{\bf y}}({\bf y}) \subset \partial \psi (S_{\epsilon})$.
Hence, thanks to \cite{FiKim} Lemma 3.1, we obtain
\begin{align}
\begin{split}\label{previousest}
\ell_d(\partial C_{{\bf y}}({\bf y}))=\ell_d\left(\partial C_{{\bf y}}({\bf y})-{\bf x}\right)\leq & C(d,\delta,R)\ell_d\left((\partial C_{{\bf y}}({\bf y})-{\bf x})\cap \delta\,\mathbb{B}_{d}\right)\\
= & C(d,\delta,R)\ell_d\left(\partial (C_{{\bf y}}({\bf y}))\cap ({\bf x}+\delta\,\mathbb{B}_d)\right)\\
\leq  & C(d,\delta,R)\ell_d\left(\partial \psi (S_{\epsilon})\cap ({\bf x}+\delta\,\mathbb{B}_d)\right)\\
\leq  & C(d,\delta,R)\ell_d\left(\partial \psi (S_{\epsilon})\cap \mathcal{X}\cap R\,\mathbb{B}_{d}\right),
\end{split}
\end{align}
where the last inequality uses the fact that 
${\bf x}+\delta'\,\mathbb{B}_{d}\subset \mathcal{X}\cap R\,\mathbb{B}_{d}$. Now we define, as before, the convex function $C^{\epsilon}_{{\bf y}}: S_{\epsilon}\rightarrow \R$ as the unique satisfying $C^{\epsilon}_{{\bf y}}({\bf y})=v_{\epsilon}(L_{\epsilon}{\bf y})$, $C_{{\bf y}}({\bf z})=\bf 0$ for ${\bf z} \in \partial L_{\epsilon}(S_{\epsilon})$, and such that ${\bf z} \mapsto C^{\epsilon}_{{\bf y}}({\bf y}-{\bf z})$ is 1-homogeneous. By the upper Alexandrov estimate (see equation (UA) in \cite{FiKim} or \cite{Fi} p.17)  applied to $v_{\epsilon}$, \begin{align*}
 |v_{\epsilon}(L_{\epsilon}({\bf y}))|^d&\leq K(d)\ell_d\left( \partial C^{\epsilon}_{{\bf y}}(L_{\epsilon}({\bf y}))\right) d_H(L_{\epsilon}(\mathbf{y}), \partial L_{\epsilon}(S_{\epsilon})),
\end{align*}
holds and hence (see eq.5 in \cite{FiKim})
\begin{align*}
|\psi_{\epsilon}({\bf y})|^d \leq K(d)\ell_d\left( \partial C_{{\bf y}}({\bf y})\right)\ell_d(S_{\epsilon})d_H(L_{\epsilon}(\mathbf{y}), \partial L_{\epsilon}(S_{\epsilon})).
\end{align*}
Then, the previous estimates \eqref{previousest} of $\ell_d(\partial C_{{\bf y}}({\bf y}))$  conclude the proof.
\end{proof}
Lemma \ref{lem:Alexandrov} and Lemma \ref{lem:bounded2}---note that we can apply the last one because we are uniformly far from $\mathbf{0}$---imply that
\begin{align*}
|v_{\epsilon}(L_{\epsilon}\bar{\bf y})|^d &\leq  \operatorname{det}(L_{\epsilon})^{{2}} C(d,\delta,R)\ell_d\left(\partial \psi (S_{\epsilon})\cap \mathcal{X}\cap R\,\mathbb{B}_{d}\right)\ell_d(S_{\epsilon})d_H(L_{\epsilon}(\mathbf{y}), \partial L_{\epsilon}(S_{\epsilon}))\\
&\leq \Lambda \operatorname{det}(L_{\epsilon})^{{2}} C(d,\delta,R)\ell_d(S_{\epsilon})^2d(L_{\epsilon}(\mathbf{y}), \partial L_{\epsilon}(S_{\epsilon}))\\
&= \Lambda  C(d,\delta,R)\ell_d(S_{\epsilon}^L)^2d(L_{\epsilon}(\mathbf{y}), \partial L_{\epsilon}(S_{\epsilon})).
\end{align*}
Since $S_{\epsilon}^L$ is normalised, 
\begin{align*}
|v_{\epsilon}(L_{\epsilon}\bar{\bf y})|^d &\leq   C(d,\delta,R, \Lambda)d_H(L_{\epsilon}(\mathbf{y}), \partial L_{\epsilon}(S_{\epsilon})),
\end{align*}
which, together with \eqref{eq:contrastFi}, implies the limit
$\lim_{\epsilon\rightarrow 0}|v_{\epsilon}(L_{\epsilon}\bar{\bf y})|=0.$ Set $A\subseteq L_{\epsilon} S_{\epsilon}$, for  $\epsilon<\epsilon_0$,   Lemma \ref{lem:previo1} implies the bound
 \begin{equation}
     \label{boundLema4.5}
     \ell_d(\partial v_{\epsilon} (A))\geq \lambda \ell_d(A) .
 \end{equation}
We consider the function 
$w({\bf z}):= \lambda^{\frac{1}{d}}\left(\frac{|{\bf z} |^2-1}{2}\right),$ 
which satisfies 
$
w=0\geq v_{\epsilon}$ in $ \partial\mathbb{B}_d, $  and $\mu_{w}(A)=\lambda\ell_d(A) \leq \ell_d(\partial v_{\epsilon} (A)),
$
where the second assertion is direct consequence of \eqref{boundLema4.5}.
Hence, Theorem 2.10. in \cite{Fi} implies that $w\geq v_{\epsilon}$ in $\mathbb{B}_d$ and consequently
$$ \left| \min _{ S_{\epsilon}^L} v_{\epsilon} \right|\geq \left| \min _{ \mathbb{B}_d} v_{\epsilon} \right| \geq \left| \min _{ \mathbb{B}_d} w\right| =\frac{\lambda^{\frac{1}{d}}}{2}=:c_0.$$
This creates a contradiction with  $\lim_{\epsilon\rightarrow 0}|v_{\epsilon}(L_{\epsilon}\bar{\bf y})|=0$ and the second limit of \eqref{eq:contrastFi}.

  \end{proof}
So far we have been dealing with the strict convexity of the function $\psi$, which implies the differentiablility of its conjugate $ \varphi $ and thus the continuity of  ${\mathbf F}_\pm$ in $\mathcal X$. 
\\ 

The continuity of the quantile map ${\mathbf Q}_\pm$ has to be restricted to the set $\Theta$. This case, as pointed out in \cite{FiKim}, can be derived by classical arguments. Indeed, the function $\psi$ is strictly convex on $\Theta$ and \emph{a fortiori} in $\Theta-\{\mathbf{0}\}$. The analogous of Corollary~4.21 in \cite{Fi}---which states that strictly convex solutions are in fact regular---can be easily derived by repeating verbatim its proof but using Lemmas~\ref{lem:Alexandrov}  and \eqref{lem:previo1} as upper and lower estimates. In this case,  for any $\mathbf{y}\in\Theta-\{\mathbf{0}\}$, while staying away from the singularity, there exist $\delta>0$ such that
\begin{align*}
	\alpha\ell_d(A)\leq \mu_{\psi}^{\mathcal{X}}(A)\leq \Lambda\ell_d(A), \ \text{for all } A\subset \mathbf{y}+\delta\mathbb{B}_d,
\end{align*}
and the same arguments of \cite{FiKim} hold. That is, strict convexity implies the continuity at $\mathbf{y}.$ Since this point was chosen arbitrarily, we can affirm that $\mathbf{Q}_{\pm}$ is continuous on $\Theta-\{\mathbf{0}\}$, which is open. To prove this last claim, let $\{\mathbf{y}_n\}_{n\in \N}\cap\Theta=\emptyset $ be a sequence  converging to $\mathbf{y}\in \mathbb{B}_d$. Lemma A.22 in \cite{Fi} yields the compactness of the set $\partial\psi(\{\mathbf{y}_n\}_{n\in \N})$. Then there exists a sub-sequence $\{\mathbf{p}_{n_k}\}_{k\in \N}$ with $\mathbf{p}_{n_k}\in \partial\psi(\mathbf{y}_{n_k})$ such that $\mathbf{p}_{n_k}\rightarrow \nabla\psi(\mathbf{y})$, which is consequence of the well known upper-semicontinuity of the subdifferential (see Exercise 12.36 in \cite{RockWets98}).  Since  $\mathbf{p}_n\not\in \mathcal{X} $, for every $n\in \N$, and $\mathcal{X}$ is open, then $\nabla\psi(\mathbf{y})\not\in \mathcal{X} $ and the claim holds. In consequence, taking smaller $\delta$ if necessary, the following  Monge-Amp\`ere inequalities hold:
\begin{align*}
	\alpha\ell_d(A)\leq \mu_{\psi}(A)\leq \Lambda\ell_d(A), \ \text{for all } A\subset \mathbf{y}+\delta\mathbb{B}_d.
\end{align*}
Therefore, standard regularity theory  proves---in the cases where $p$ is regular enough---the continuity of the derivatives. Moreover,  the map $\mathbf{F}_{\pm}$ defines an  homeomorphism between $\mathcal{X}-\{\mathbf{Q}_{\pm}(\mathbf{0})\}$ and $\Theta-\{\mathbf{0}\}$, whose inverse is the quantile $\mathbf{Q}_{\pm}$. 

\qed
\section*{Acknowledgments}
The authors would like to thank Johan Segers for his insightful discussions and for providing reference \cite{Segers} for this work.
\section*{Funding}
The research of Eustasio del Barrio and Alberto González Sanz is partially supported by grant PID2021-128314NB-I00 funded by
MCIN/AEI/ 10.13039/501100011033/FEDER, UE. The research of Alberto Gonz\'alez Sanz is partially supported by the AI Interdisciplinary Institute ANITI, which is funded by the French \textit{``investing for the Future – PIA3''} program under the Grant agreement ANR-19-PI3A-0004. 

\bibliographystyle{elsarticle-harv} 
\bibliography{ref}


\end{document}